\documentclass[a4paper,draft,leqno,11pt]{article}
\usepackage{ifthen,amsfonts,amssymb,epic,eepic,stmaryrd,mathrsfs,bm}
\usepackage{amsmath,a4wide}
\usepackage{theorem}
\usepackage[pdftex]{graphics,color}

\makeatletter\@addtoreset{equation}{section}\makeatother

\newtheorem{theorem}{Theorem}[section]

\begingroup
\newtheorem{lemma}{Lemma}[section]
\newtheorem{proposition}[lemma]{Proposition}

\endgroup

\newtheorem{definition}{Definition}[section]
{\theorembodyfont{\rmfamily}}

\newcommand{\ep}{\varepsilon}

\newcommand{\ds}{\displaystyle}
\newcommand{\beq}[1]{\begin{equation} \label{#1}\ds}
\newcommand{\eeq}{\end{equation}}
\newcommand{\bml}[1]{\beq{#1} \begin{array}{c}\ds}
\newcommand{\eml}{\end{array}\eeq}
\newcommand{\beqq}{\begin{equation*}\ds}
\newcommand{\eeqq}{\end{equation*}}
\newcommand{\bmll}{\beqq \begin{array}{c}\ds}
\newcommand{\emll}{\end{array}\eeqq}
\renewcommand{\div}{{\rm div}\,}

\newcommand{\abs}[1]{\ensuremath{\left| #1 \right|}}

\def \de{\partial}

\def \ep{\varepsilon}

\def \d{\mathrm{d}}

\newcommand{\R}{\mathbb{R}}

\newcommand{\id}{{\rm Id}}
\newcommand{\Sym}{\mathcal{S}}

\begin{document}

%opening
\author{Elisabetta Chiodaroli and Ond\v{r}ej Kreml}
\title{An overview of some recent results on the Euler system of isentropic gas dynamics }
\date{}

\maketitle

\centerline{EPFL Lausanne}

\centerline{Station 8, CH-1015 Lausanne, Switzerland}

\bigskip

\centerline{Institute of Mathematics of the Academy of Sciences of the Czech Republic}

\centerline{\v{Z}itn\'a 25, 115 67 Praha 1, Czech Republic}

%\maketitle
 
\begin{abstract}

This overview is concerned with the well-posedness problem for the isentropic
compressible Euler equations of gas dynamics. The results we
present are in line with the program of investigating the efficiency of different 
selection criteria proposed in the literature in order to weed out non-physical
solutions to more-dimensional systems of conservation laws and they build upon
the method of convex integration developed by De Lellis and Sz\'ekelyhidi for the 
incompressible Euler equations. 
Mainly following \cite{chk}, we investigate the role of the maximal dissipation criterion proposed by Dafermos in \cite{Da1}:
we prove how, for specific pressure laws, some non-standard (i.e. constructed via convex integration methods) solutions to the Riemann problem for the isentropic
Euler system in two space dimensions have greater energy dissipation rate than
the classical self-similar solution emanating from the same Riemann data. We
therefore show that the maximal dissipation criterion proposed by Dafermos
does not favour in general the self-similar solutions.
\end{abstract}
{\let\thefootnote\relax

\section{Introduction}

We consider the isentropic compressible Euler
system of gas dynamics in $2$ space dimensions (cf. \cite{da} or \cite{se} or \cite{br}). It is obtained as a simplification
of the full compressible Euler equations, by assuming the entropy to be constant.
The state of the gas is described through the state vector
$$ V=(\rho, v)$$
whose components are the density $\rho$ and the velocity $v$. The system consists of $3$ equations which correspond to balance statements for mass and linear momentum.
The corresponding Cauchy problems reads as
\begin{equation}\label{eq:Euler system}
\left\{\begin{array}{l}
\partial_t \rho + {\rm div}_x (\rho v) \;=\; 0\\
\partial_t (\rho v) + {\rm div}_x \left(\rho v\otimes v \right) + \nabla_x [ p(\rho)]\;=\; 0\\
\rho (\cdot,0)\;=\; \rho^0\\
v (\cdot, 0)\;=\; v^0 \, ,
\end{array}\right.
\end{equation}
with $t\in \R^+$, $x\in \R^2$.
The pressure $p$ is a function of the density $\rho$ determined from the constitutive thermodynamic
relations of the gas under consideration and it is assumed to satisfy $p'>0$ (under this assumption the system is hyperbolic).
We will work with pressure laws
$p(\rho)=  \rho^\gamma$
with constant $\gamma\geq 1$. 

Our aim is to discuss the issue of uniqueness of weak solutions to the Cauchy problem \eqref{eq:Euler system}.
The theory of the Cauchy problem for hyperbolic systems of conservation laws is typically confronted with two major challenges.
First, it is well-known that classical solutions develop discontinuities, even starting out from smooth initial data. 
In the literature this behaviour is known as breakdown of classical solutions.
Therefore, it becomes imperative to introduce the notion of weak solution. However, weak solutions fail to be unique. 
In order to restore uniqueness restrictions need to be imposed in hope of singling
out a unique physical solution. When dealing with systems coming from Physics, as in our case, the second law of Thermodynamics naturally induces such restrictions, 
such admissibility criteria by stipulating that weak solutions are admissible/entropy solutions if they satisfy some \textit{entropy} inequalities 
(see \eqref{eq:energy inequality} for the specific case of the compressible Euler system).
Finally, a third important challenge then arises: do \textit{entropy} inequalities really serve as selection criteria? Are admissible solutions unique? Or at least, do there exist 
efficient criteria restoring uniqueness?
This is a central problem to set down a complete theory for the Cauchy problem. It has deserved a lot of attention, but positive answers were found only
for scalar conservation laws or for systems in one space dimensions (under smallness assumptions on the initial data): for a complete account of the existing literature we refer the reader 
to \cite{da} and \cite{se}.
When dealing with systems of conservation laws in more than one space dimension, it is still an
intriguing mathematical problem to develop a theory of well-posedness for
the Cauchy problem which includes the formation and evolution of
shock waves. 

In 2006, Elling \cite{el} studied numerically a particular case of initial data for the two dimensional non-isentropic Euler
equations. His results show that the numerical method
does not always converge to the physical solution.
Moreover, they suggest that entropy solutions (in the weak entropy inequality
sense) to the multi-dimensional Euler equations are not always unique.

In a groundbreaking paper \cite{dls2}, De Lellis-Sz\'{e}kelyhidi give an example in
favour of the conjecture that entropy/admissible solutions to the
multi-dimensional compressible Euler equations are in general not
unique. The non-uniqueness result by De Lellis-Sz\'{e}kelyhidi is a byproduct of
their analysis of the incompressible Euler equations based on
its formulation as a differential inclusion (see \cite{dls1} and \cite{dls3}) combined with convex integration methods: they exploit the result for the incompressible Euler
equations to exhibit bounded initial density and bounded compactly supported initial velocity for which admissible
solutions of \eqref{eq:Euler system} are not unique (in more than one space dimension).
However the initial data constructed in \cite{dls2} are very irregular. The result by De Lellis and Sz\'ekelyhidi is improved by the author in \cite{ch}
 where it is proven that non-uniqueness still holds in the case of regular initial density (see also \cite{ChFeKr} for further generalizations). Non--unique solutions constructed via convex integration are 
referred to as non--standard or oscillatory solutions.
Moreover in \cite{ChDLKr}, using the Riemann problem as a building block, the authors show that, in the two dimensional case, the entropy inequality (see
\eqref{eq:energy inequality})
does not single out unique weak solutions even under very strong assumptions on the initial data ($(\rho^0, v^0) \in W^{1, \infty}(\R^2)$):
\begin{theorem}[Chiodaroli, De Lellis, Kreml]\label{t:lipschitz}
There are Lipschitz initial data $(\rho^0, v^0)$ for which there are infinitely
many bounded admissible solutions $(\rho, v)$ of \eqref{eq:Euler system} on
$\R^2\times [0, \infty[$ with
$\inf \rho >0$. These
solutions are all locally Lipschitz on a finite interval on which they
all coincide with the unique classical solution.   
\end{theorem}
This is proven by constructing infinitely many entropy weak solutions in forward time to a Riemann problem for \eqref{eq:Euler system}
whose Riemann data can be generated, backwards in time, by a classical compression wave: the Lipschitz initial data of Theorem \ref{t:lipschitz}
will be provided by the values of the compression wave at some finite negative time. It is clear now that the infinitely many admissible solutions
constructed in Theorem \ref{t:lipschitz} all coincide with the unique classical solution (compression wave) on a finite time interval
 whereas non--uniqueness arises after the first blow--up time.

This series of negative results concerning the entropy inequality as selection criterion for system \eqref{eq:Euler system} motivated 
the authors to explore other admissibility criteria which could work in favour of uniqueness,
in particular we investigated an alternative criterion which has been proposed by Dafermos in \cite{Da1} under
the name \textit{entropy rate admissibility criterion}.
The ideas developed in \cite{ChDLKr} enabled us to prove in \cite{chk} the following theorem:
\begin{theorem}[Chiodaroli, Kreml]\label{t:mainthm}
 Let $p(\rho) = \rho^\gamma$, $1 \leq \gamma < 3$. 
 There exist Riemann data for which the self-similar solution to \eqref{eq:Euler system}
 emanating from these data is not entropy rate admissible.
\end{theorem}
This result does not exclude that the entropy rate admissibility criterion could still select a unique solution, but surely prevents the self--similar solution
to be the selected one. Moreover, since Theorem \ref{t:mainthm} is proven using non--standard solutions as competitors, with respect to Dafermos' criterion,
for the self--similar solution, we can affirm that the entropy rate criterion cannot, at least in our setting, rule out oscillatory solutions obtained via convex integration.
 In the rest of the paper we will further explain this result.

\section{Entropy criteria}

\subsection{Entropy inequality}
We recall here the usual definitions of weak and admissible
solutions to \eqref{eq:Euler system} in the two--dimensional case.

\begin{definition}\label{d:weak}
By a \emph{weak solution} of \eqref{eq:Euler system} on $\R^2\times[0,\infty[$ we
mean a pair $(\rho, v)\in L^\infty(\R^2\times [0,\infty[)$ such that the following identities 
hold for every test functions $\psi\in C_c^{\infty}(\R^2\times [0, \infty[)$,
$\phi\in C_c^{\infty}(\R^2\times [0, \infty[)$:
\begin{equation} \label{eq:weak1}
\int_0^\infty \int_{\R^2} \left[\rho\partial_t \psi+ \rho v \cdot \nabla_x \psi\right] dx dt+\int_{\R^2} \rho^0(x)\psi(x,0) dx \;=\; 0
\end{equation}
\begin{align} \label{eq:weak2}
&\int_0^\infty \int_{\R^2} \left[ \rho v \cdot \partial_t \phi+ \rho v \otimes v : D_x \phi +p(\rho) \div_x \phi \right]
+\int_{\R^2} \rho^0(x) v^0(x)\cdot\phi(x,0) dx\;=\; 0.
\end{align}
\end{definition}

\begin{definition}\label{d:admissible}
A bounded weak solution $(\rho, v)$ of \eqref{eq:Euler system} is \emph{admissible} if
it satisfies the following inequality for every nonnegative 
test function $\varphi\in C_c^{\infty}(\R^2\times [0,\infty[)$:
\begin{align} \label{eq:energy inequality}
 &\int_0^\infty\int_{\R^2} \left[\left(\rho\varepsilon(\rho)+\rho \frac{\abs{v}^2}{2}\right)\partial_t \varphi+\left(\rho\varepsilon(\rho)+\rho
\frac{\abs{v}^2}{2}+p(\rho)\right) v \cdot \nabla_x \varphi \right]\notag \\
&+\int_{\R^2} \left(\rho^0 (x) \varepsilon(\rho^0 (x))+\rho^0 (x)\frac{\abs{v^0 (x)}^2}{2}\right) 
\varphi(x,0)\, dx \;\geq\; 0\, .
\end{align}
\end{definition}
Note that \eqref{eq:energy inequality} is rather a weak form of energy balance. 

% \begin{figure}[htbp]
% \begin{center}
% \input fig.pdf_t
% \caption{shock}
% \label{fig:fan}
% \end{center}
% \end{figure}

\subsection{Entropy rate admissibility criterion} \label{s:entropy}

An alternative criterion to the entropy inequality has been proposed by Dafermos in \cite{Da1} under
the name of \textit{entropy rate admissibility criterion}.
In order to formulate this criterion for the specific system \eqref{eq:Euler system} we define the \textit{total energy} 
of the solutions $(\rho,v)$ to \eqref{eq:Euler system} as
\begin{equation} \label{eq:energy0}
E[\rho,v](t) = \int_{\R^2} \left(\rho\ep(\rho) + \rho\frac{\abs{v}^2}{2}\right)\d x .
\end{equation}
Let us remark that in Dafermos' terminology $E[\rho,v](t)$ is called ``total entropy'' (see \cite{Da1}). However, since in the context of system \eqref{eq:Euler system}
the physical energy plays the role of the mathematical entropy, it is more convenient to call $E[\rho,v](t)$ \textit{total energy}.
The right derivative of $E[\rho,v](t)$ defines the \textit{energy dissipation rate} of $(\rho,v)$ at time $t$:
\begin{equation} \label{eq:dissipation rate0}
D[\rho,v](t) = \frac{\d_+ E[\rho,v](t)}{\d t}. 
\end{equation}
Since our solutions will have piecewise constant values of $\rho$ and $\abs{v}^2$ and it is easy to see that the total energy of any solution 
we construct is infinite, we shall restrict the infinite domain $\R^2$ to a finite box $(-L,L)^2$ and denote
\begin{align}
&E_L[\rho,v](t) = \int_{(-L,L)^2} \left(\rho\ep(\rho) + \rho\frac{\abs{v}^2}{2}\right)\d x \label{eq:energy L}\\
&D_L[\rho,v](t) = \frac{\d_+ E_L[\rho,v](t)}{\d t}. \label{eq:dissipation rate L}
\end{align}
The problem of infinite energy of solutions may be solved also by restricting to a periodic domain
and constructing (locally in time) periodic solutions. This procedure is carefully described in \cite{chk}.

Using the concept of energy dissipation rate, Dafermos in \cite{Da1} introduces a new selection criterion for weak solutions
which goes under the name of \textit{entropy rate admissibility criterion}. We recall here the definition of \textit{entropy rate admissible solutions}.
\begin{definition}[Entropy rate admissible solution]\label{d:entropy rate}
A weak solution $(\rho,v)$ of \eqref{eq:Euler system} is 
\emph{entropy rate admissible} if 
there exists $L^* > 0$ such that there is no other weak solution $(\overline{\rho},\overline{v})$ with the property that for some $\tau\geq 0$, $(\overline{\rho},\overline{v})(x,t)= (\rho,v)(x,t)$ on $\R^2 \times [0, \tau]$
and $ D_L[\overline{\rho},\overline{v}](\tau) < D_L[\rho,v](\tau) $ for all $L \geq L^*$.
\end{definition}
In other words, we call entropy rate admissible the solution(s) dissipating most total energy.

\section{Background literature and main results}
The investigation of the entropy rate admissibility criterion initiated with the paper \cite{Da1}
of Dafermos where he puts it forward and moreover proves that for a single equation the entropy rate criterion is equivalent to the viscosity criterion in the class
of piecewise smooth solutions. 
Later on, following the approach of Dafermos, Hsiao in \cite{Hs} proves, in the class
of piecewise smooth solutions, the equivalence of the entropy rate criterion and
the viscosity criterion for the one-dimensional system of equations of nonisentropic gas dynamics in lagrangian formulation with pressure laws $p(\rho)= \rho^\gamma$
for $\gamma\geq 5/3 $ while the same equivalence is disproved for $\gamma < 5/3 $.
For further analysis on the relation between entropy rate minimization and admissibility of solutions for a more general class of evolutionary equations we refer to \cite{Da2}.
However, to our knowledge, up to some time ago the entropy rate criterion had not been tested in the case of several space variables and on broader class of solutions
than the piecewise smooth ones.

Recently Feireisl in \cite{fe} extended the result of Chiodaroli \cite{ch} and obtained infinitely many global admissible weak solutions of \eqref{eq:Euler system}
none of which is entropy rate admissible: this results seems in favour of the effectiveness of the entropy rate criterion 
to rule out non--standard solutions (i.e. constructed by the method of De Lellis and Sz\'ekelyhidi).
In \cite{chk} we have actually shown that for specific initial data, and in the two--dimensional case,
the oscillatory (non--standard) solutions dissipate more energy than the self-similar solution which may be believed to be the physical one.
The results obtained in \cite{ChDLKr} hinge upon some of the ideas devised in \cite{ChDLKr} combined with novel developments to deal with
the entropy rate criterion.

We refer also to the work \cite{sz}, where Sz\'ekelyhidi
constructed irregular solutions of the incompressible Euler equations with vortex-sheet initial data and computed their dissipation rate.

We focus on the Riemann problem for the system \eqref{eq:Euler system},\eqref{eq:energy inequality} in two-space dimensions. 
Hence, we denote the space variable as $x=(x_1, x_2)\in \R^2$ and consider initial data in the form
\begin{equation}\label{eq:R_data}
(\rho^0 (x), v^0 (x)) := \left\{
\begin{array}{ll}
(\rho_-, v_-) \quad & \mbox{if $x_2<0$}\\ \\
(\rho_+, v_+) & \mbox{if $x_2>0$,} 
\end{array}\right. 
\end{equation}
where $\rho_\pm, v_{\pm 1}, v_{\pm 2}$ are constants. 
Our concern has been to compare the energy dissipation rate of standard self-similar solutions associated to the Riemann problem \eqref{eq:Euler system}, \eqref{eq:energy inequality},
\eqref{eq:R_data} with the energy dissipation rate of non-standard solutions for the same problem obtained by the method developed in \cite{ChDLKr}. 

We obtained the following results.

\begin{theorem}\label{t:main0}
 Let $p(\rho) = \rho^{\gamma}$ with $\gamma \geq 1$. For every Riemann data \eqref{eq:R_data} such that the self-similar solution 
 to the Riemann problem \eqref{eq:Euler system}--\eqref{eq:energy inequality}, \eqref{eq:R_data} consists of an admissible $1-$shock and an admissible $3-$shock, i.e. $v_{-1} = v_{+1}$ and
 \begin{equation}\label{eq:2shocks condition}
  v_{+2} - v_{-2} < -\sqrt{\frac{(\rho_--\rho_+)(p(\rho_-)-p(\rho_+))}{\rho_-\rho_+}},
 \end{equation}
 there exist infinitely many admissible solutions to \eqref{eq:Euler system}--\eqref{eq:energy inequality}, \eqref{eq:R_data}.
\end{theorem}

Theorem \ref{t:main0} can be viewed as an extension of the results obtained 
together with De Lellis in \cite{ChDLKr}. 
As a consequence of Theorem \ref{t:main0} and by a suitable choice of initial data, we can prove the following main theorem.

\begin{theorem}\label{t:main}
 Let $p(\rho) = \rho^\gamma$, $1 \leq \gamma < 3$. 
 There exist Riemann data \eqref{eq:R_data} for which the self-similar solution to \eqref{eq:Euler system},\eqref{eq:energy inequality}
 emanating from these data is not entropy rate admissible.
\end{theorem}

Theorem \ref{t:main} ensures that for $1 \leq \gamma < 3$ there exist initial Riemann data \eqref{eq:R_data} for which 
some of the infinitely many nonstandard solutions constructed as in Theorem \ref{t:main0} dissipate more energy than the self-similar solution, 
suggesting in particular that the Dafermos entropy rate admissibility criterion would not pick the self-similar solution as the admissible one.

\section{Strategies of proof}\label{s:I}

In this section we explain how to prove Theorem \ref{t:main0} and \ref{t:main}.
For the complete proofs we refer the reader to \cite{chk}.

Both theorems stem from the framework introduced in \cite{ChDLKr} where the Riemann problem constitutes the starting point for constructing
non--unique admissible non--standard solutions. In particular, in \cite{ChDLKr}, the authors jointly with Camillo De Lellis developed
a method which allows to obtain infinitely many entropy (oscillatory) solutions to a Riemann problem provided a suitable \textit{admissible subsolution}
exists.

\subsection{From subsolutions to solutions}

We shall introduce the notion of \textit{fan subsolution}  and \textit{admissible fan subsolution} as in \cite[Section 3]{ChDLKr}.

\begin{definition}[Fan partition]\label{d:fan}
A {\em fan partition} of $\R^2\times (0, \infty)$ consists of three open sets $P_-, P_1, P_+$
of the following form 
\begin{align}
 P_- &= \{(x,t): t>0 \quad \mbox{and} \quad x_2 < \nu_- t\}\\
 P_1 &= \{(x,t): t>0 \quad \mbox{and} \quad \nu_- t < x_2 < \nu_+ t\}\\
 P_+ &= \{(x,t): t>0 \quad \mbox{and} \quad x_2 > \nu_+ t\},
\end{align}
where $\nu_- < \nu_+$ is an arbitrary couple of real numbers.
\end{definition}

We denote by $\Sym_0^{2\times2}$ the set of all symmetric $2\times2$ matrices with zero trace.

\begin{definition}[Fan subsolution] \label{d:subs}
A {\em fan subsolution} to the compressible Euler equations \eqref{eq:Euler system} with
initial data \eqref{eq:R_data} is a triple 
$(\overline{\rho}, \overline{v}, \overline{u}): \R^2\times 
(0,\infty) \rightarrow (\R^+, \R^2, \Sym_0^{2\times2})$ of piecewise constant functions satisfying
the following requirements.
\begin{itemize}
\item[(i)] There is a fan partition $P_-, P_1, P_+$ of $\R^2\times (0, \infty)$ such that
\[
(\overline{\rho}, \overline{v}, \overline{u})=  
(\rho_-, v_-, u_-) \bm{1}_{P_-}
+ (\rho_1, v_1, u_1) \bm{1}_{P_1}
+ (\rho_+, v_+, u_+) \bm{1}_{P_+}
\]
where $\rho_1, v_1, u_1$ are constants with $\rho_1 >0$ and $u_\pm =
v_\pm\otimes v_\pm - \textstyle{\frac{1}{2}} |v_\pm|^2 \id$;
\item[(ii)] There exists a positive constant $C$ such that
\begin{equation} \label{eq:subsolution 2}
v_1\otimes v_1 - u_1 < \frac{C}{2} \id\, ;
\end{equation}
\item[(iii)] The triple $(\overline{\rho}, \overline{v}, \overline{u})$ solves the following system in the
sense of distributions:
\begin{align}
&\partial_t \overline{\rho} + {\rm div}_x (\overline{\rho} \, \overline{v}) \;=\; 0\label{eq:continuity}\\
&\partial_t (\overline{\rho} \, \overline{v})+{\rm div}_x \left(\overline{\rho} \, \overline{u} 
\right) + \nabla_x \left( p(\overline{\rho})+\frac{1}{2}\left( C \rho_1
\bm{1}_{P_1} + \overline{\rho} |\overline{v}|^2 \bm{1}_{P_+\cup P_-}\right)\right)= 0.\label{eq:momentum}
\end{align}
\end{itemize}
\end{definition}

\begin{definition}[Admissible fan subsolution]\label{d:admiss}
 A fan subsolution $(\overline{\rho}, \overline{v}, \overline{u})$ is said to be {\em admissible}
if it satisfies the following inequality in the sense of distributions
\begin{align} 
&\de_t \left(\overline{\rho} \varepsilon(\overline{\rho})\right)+\div_x
\left[\left(\overline{\rho}\varepsilon(\overline{\rho})+p(\overline{\rho})\right) \overline{v}\right]
 + \de_t \left( \overline{\rho} \frac{|\overline{v}|^2}{2} \bm{1}_{P_+\cup P_-} \right)
+ \div_x \left(\overline{\rho} \frac{|\overline{v}|^2}{2} \overline{v} \bm{1}_{P_+\cup P_-}\right)\nonumber\\
&\qquad\qquad+ \left[\de_t\left(\rho_1 \, \frac{C}{2} \, \bm{1}_{P_1}\right) 
+ \div_x\left(\rho_1 \, \overline{v} \, \frac{C}{2}  \, \bm{1}_{P_1}\right)\right]
\;\leq\; 0\, .\label{eq:admissible subsolution}
\end{align}
\end{definition}

The strategy which lies behind Theorem \ref{t:main0}, as well as behind Theorem \ref{t:lipschitz} in \cite{ChDLKr},
consists in reducing the proof of the existence of infinitely many admissible solutions to the Riemann problem 
for \eqref{eq:Euler system} to the proof of 
the existence of an admissible fan subsolution as defined in Definition \ref{d:admiss}.
This is the content of the following Proposition
which represents the key ingredient of \cite{chk} and is proven in \cite{ChDLKr}. 

\begin{proposition}\label{p:subs}
Let $p$ be any $C^1$ function and $(\rho_\pm, v_\pm)$ be such that there exists at least one
admissible fan subsolution $(\overline{\rho}, \overline{v}, \overline{u})$ of \eqref{eq:Euler system}
with initial data \eqref{eq:R_data}. Then there are infinitely 
many bounded admissible solutions $(\rho, v)$ to \eqref{eq:Euler system}-\eqref{eq:energy inequality}, \eqref{eq:R_data} such that 
$\rho=\overline{\rho}$ and $\abs{v}^2\bm{1}_{P_1} = C$.
\end{proposition}

Roughly speaking, the infinitely many bounded admissible solutions $(\rho,v)$ of Proposition \ref{p:subs}
are constructed by adding to the subsolution solutions to 
the linearized pressureless incompressible Euler equations supported in $P_1$. Such solutions are given by the following Lemma, cf. \cite[Lemma 3.7]{ChDLKr}.

\begin{lemma}\label{l:ci}
Let $(\tilde{v}, \tilde{u})\in \R^2\times \Sym_0^{2\times 2}$ and $C_0>0$ be such that $\tilde{v}\otimes \tilde{v}
- \tilde{u} < \frac{C_0}{2} \id$. For any open set $\Omega\subset \R^2\times \R$ there are infinitely many maps
$(\underline{v}, \underline{u}) \in L^\infty (\R^2\times \R , \R^2\times \Sym_0^{2\times 2})$ with the following property
\begin{itemize}
\item[(i)] $\underline{v}$ and $\underline{u}$ vanish identically outside $\Omega$;
\item[(ii)] $\div_x \underline{v} = 0$ and $\partial_t \underline{v} + \div_x \underline{u} = 0$;
\item[(iii)] $ (\tilde{v} + \underline{v})\otimes (\tilde{v} + \underline{v}) - (\tilde{u} + \underline{u}) = \frac{C_0}{2} \id$
a.e. on $\Omega$.
\end{itemize}
\end{lemma}

Proposition \ref{p:subs} is then proved by applying Lemma \ref{l:ci} with $\Omega = P_1$, $(\tilde{v}, \tilde{u}) = (v_1,u_1)$ and $C_0 = C$. It is then a matter of easy computation to check that each couple $(\overline{\rho}, \overline{v} + \underline{v})$ is indeed an admissible weak solution to \eqref{eq:Euler system}--\eqref{eq:energy inequality} with initial data \eqref{eq:R_data}, for details see \cite[Section 3.3]{ChDLKr}.
The whole proof of Lemma \ref{l:ci} can be found in \cite[Section 4]{ChDLKr}.

\subsection{Concluding arguments}\label{s:concluding}

Thanks to Proposition \ref{p:subs}, Theorem \ref{t:main0} amounts to showing the existence of a fan admissible subsolution with appropriate initial data
under the hypothesis that \eqref{eq:R_data} is such that the self-similar solution 
 to the Riemann problem \eqref{eq:Euler system}, \eqref{eq:energy inequality}, \eqref{eq:R_data} consists of an admissible $1-$shock and an admissible $3-$shock

Indeed, a fan admissible sunsolution with initial data \eqref{eq:R_data}
is defined through a the set of identities and inequalities which we recall here (see also \cite[Section 5]{ChDLKr}).

We introduce the real numbers 
$\alpha, \beta, \gamma_1, \gamma_2, v_{-1}, v_{-2}, v_{+1}, v_{+2}$ such that
\begin{align} 
v_1 &= (\alpha, \beta),\label{eq:v1}\\
v_- &= (v_{-1}, v_{-2})\\
v_+ &= (v_{+1}, v_{+2})\\
u_1 &=\left( \begin{array}{cc}
    \gamma_1 & \gamma_2 \\
    \gamma_2 & -\gamma_1\\
    \end{array} \right)\, .\label{eq:u1}
\end{align}

Then, Proposition \ref{p:subs} translates into the following set of algebraic identities and inequalities.

\begin{proposition}\label{p:algebra}
Let $P_-, P_1, P_+$ be a fan partition as in Definition \ref{d:fan}.\\
 The constants 
$v_1, v_-, v_+, u_1, \rho_-, \rho_+, \rho_1$ as in \eqref{eq:v1}-\eqref{eq:u1} define an 
\emph{admissible fan subsolution} as in Definitions \ref{d:subs}-\ref{d:admiss} if and only if the following
identities and inequalities hold:
\begin{itemize}
\item Rankine-Hugoniot conditions on the left interface:
\begin{align}
&\nu_- (\rho_- - \rho_1) \, =\,  \rho_- v_{-2} -\rho_1  \beta \label{eq:cont_left}  \\
&\nu_- (\rho_- v_{-1}- \rho_1 \alpha) \, = \, \rho_- v_{-1} v_{-2}- \rho_1 \gamma_2  \label{eq:mom_1_left}\\
&\nu_- (\rho_- v_{-2}- \rho_1 \beta) \, = \,  
\rho_- v_{-2}^2 + \rho_1 \gamma_1 +p (\rho_-)-p (\rho_1) - \rho_1 \frac{C}{2}\, ;\label{eq:mom_2_left}
\end{align}
\item Rankine-Hugoniot conditions on the right interface:
\begin{align}
&\nu_+ (\rho_1-\rho_+ ) \, =\,  \rho_1  \beta - \rho_+ v_{+2} \label{eq:cont_right}\\
&\nu_+ (\rho_1 \alpha- \rho_+ v_{+1}) \, = \, \rho_1 \gamma_2 - \rho_+ v_{+1} v_{+2} \label{eq:mom_1_right}\\
&\nu_+ (\rho_1 \beta- \rho_+ v_{+2}) \, = \, - \rho_1 \gamma_1 - \rho_+ v_{+2}^2 +p (\rho_1) -p (\rho_+) 
+ \rho_1 \frac{C}{2}\, ;\label{eq:mom_2_right}
\end{align}
\item Subsolution condition:
\begin{align}
 &\alpha^2 +\beta^2 < C \label{eq:sub_trace}\\
& \left( \frac{C}{2} -{\alpha}^2 +\gamma_1 \right) \left( \frac{C}{2} -{\beta}^2 -\gamma_1 \right) - 
\left( \gamma_2 - \alpha \beta \right)^2 >0\, ;\label{eq:sub_det}
\end{align}
\item Admissibility condition on the left interface:
\begin{align}
& \nu_-(\rho_- \varepsilon(\rho_-)- \rho_1 \varepsilon( \rho_1))+\nu_- 
\left(\rho_- \frac{\abs{v_-}^2}{2}- \rho_1 \frac{C}{2}\right)\nonumber\\
\leq & \left[(\rho_- \varepsilon(\rho_-)+ p(\rho_-)) v_{-2}- 
( \rho_1 \varepsilon( \rho_1)+ p(\rho_1)) \beta \right] 
+ \left( \rho_- v_{-2} \frac{\abs{v_-}^2}{2}- \rho_1 \beta \frac{C}{2}\right)\, ;\label{eq:E_left}
\end{align}
\item Admissibility condition on the right interface:
\begin{align}
&\nu_+(\rho_1 \varepsilon( \rho_1)- \rho_+ \varepsilon(\rho_+))+\nu_+ 
\left( \rho_1 \frac{C}{2}- \rho_+ \frac{\abs{v_+}^2}{2}\right)\nonumber\\
\leq &\left[ ( \rho_1 \varepsilon( \rho_1)+ p(\rho_1)) \beta- (\rho_+ \varepsilon(\rho_+)+ p(\rho_+)) v_{+2}\right] 
+ \left( \rho_1 \beta \frac{C}{2}- \rho_+ v_{+2} \frac{\abs{v_+}^2}{2}\right)\, .\label{eq:E_right}
\end{align}
\end{itemize}
\end{proposition}

Theorem \ref{t:main} is then a corollary of the following theorem, proven in \cite{chk} and showing the existence of an admissible fan subsolution, combined with Proposition \ref{p:subs}.

\begin{theorem} \label{t:intermediate}
  Let $p(\rho) = \rho^{\gamma}$ with $\gamma \geq 1$.
 For every Riemann data \eqref{eq:R_data} such that $v_{-1} = v_{+1}$ and
 \begin{equation}\label{eq:2shocks condition}
  v_{+2} - v_{-2} < -\sqrt{\frac{(\rho_--\rho_+)(p(\rho_-)-p(\rho_+))}{\rho_-\rho_+}},
 \end{equation}
 there exist $\nu_{-}, \nu_+, v_1, u_1,\rho_1, C$ such that \eqref{eq:cont_left}--\eqref{eq:E_right} hold.
\end{theorem}

It remains to prove Theorem \ref{t:main}. This is obtained by showing that 
among the infinitely many admissible solutions provided by Theorem 
\ref{t:main0} one has lower energy dissipation rate than the self--similar solution emanating from the same Riemann data, thus contradicting Definition \ref{d:entropy rate}.
Let us remark that the Riemann data allowing for the result of Theorem \ref{t:main} are of the same type as the ones of Theorem \ref{t:main0}, i.e.
they admit a forward in time self--similar solution consisting of two shocks. We also underline that the self--similar solution depends in fact only on one variable,
specifically on $x_2$.

Assume from now on for simplicity that $ v_{-1} = v_{+1} = 0$ in \eqref{eq:R_data}.
Let us denote the self--similar solution emanating from the Riemann data as in Theorem \ref{t:main0}.
The value of the dissipation rate $D_L[\rho_S,v_S](t)$ has a specific form for the solution $(\rho_S, v_S)$ 
consisting (by assumption) of two shocks of speeds $\nu_1$ and $\nu_2$. Denoting the middle state $(\rho_m,v_m = (0,\overline{v}))$ and introducing the notation
\begin{align}
 E_\pm &:= \rho_\pm\ep(\rho_\pm) + \rho_\pm\frac{v_\pm^2}{2} \label{eq:E pm} \\
 E_m &:= \rho_m\ep(\rho_m) + \rho_m\frac{\overline{v}^2}{2}
\end{align}
we have
\begin{equation}\label{eq:Dis rate for 2 shocks}
 D_L[\rho_S,v_S] = -2L\left(\nu_1(E_- -E_m) + \nu_2(E_m-E_+)\right).
\end{equation}
Now let us consider a solution $(\rho, v)$ with the same initial data \eqref{eq:R_data} but constructed by the method of convex integration starting 
from an admissible fan subsolution using Proposition \ref{p:subs}. We also assume, that the fan admissible subsolution (which exists by 
Theorem \ref{t:intermediate}) has an underlying fan partition defined by the speeds $\nu_-$ and $\nu_+$.
Although $v$ is not constant in $P_1$ we still have, by construction (see Proposition \ref{p:subs}), 
that $\abs{v}^2\bm{1}_{P_1} = C $, in particular $E_1: = \rho_1\ep(\rho_1) + \rho_1\frac{C}{2}$ 
is constant in $P_1$. The dissipation rate for all solutions constructed from a given subsolution hence
depends only on the underlying subsolution and is given by 
\begin{equation}\label{eq:Dis rate for non-standard}
 D_L[\rho,v] = -2L\left(\nu_-(E_- -E_1) + \nu_+(E_1-E_+)\right).
\end{equation} 
If, for a moment, we assume that the speeds of the self--similar solution and of the subsolution coincide, i.e. 
$\nu_-=\nu_1$ and $\nu_+=\nu_2$, it is clearly enough to achieve $E_1>E_m$ in order to prove Theorem \ref{t:main}.
Of course, as one can see from \cite[Section 4]{chk} this is not the case; nonetheless the proof works
along the same line: we can prove that there is still some freedom in the choice
of the parameters defining the underlying subsolution for $(\rho,v)$ which allows to tune them in such a way that Theorem \ref{t:main} holds.
For a complete and rigorous proof we refer to \cite[Section 5]{chk}.

\end{document}